\documentclass[a4paper]{article}

\usepackage{a4}
\usepackage{amsmath, amssymb}   %
\usepackage{mathtools}
\usepackage{graphicx,color}                   %
\usepackage{float}  
\usepackage{xspace}         
\usepackage[scientific-notation=true]{siunitx}         

\usepackage{tikz}
\tikzstyle{startstop} = [rectangle, rounded corners,minimum width=3cm, minimum height=1cm, text centered, text width=5cm, draw=black]
\tikzstyle{process} = [rectangle, minimum width=3cm, minimum height=1cm, text width=6 cm,text centered, draw=black]
\tikzstyle{decision} = [diamond, minimum width=3cm, minimum height=1cm, text width=3cm, text centered, draw=black]
\tikzstyle{arrow} = [thick,->,>=stealth]
\usetikzlibrary{spy}
\usetikzlibrary{shapes.geometric, arrows}
\usetikzlibrary{positioning} 
\usetikzlibrary{calc}
\usepackage{pgfplots}
\pgfplotsset{width=10cm,compat=1.9}
\usepackage[a4paper, total={6.5in, 8.5in}]{geometry}

\usepackage{hyperref}
\usepackage[all]{hypcap}

\definecolor{light_gray}{gray}{0.75}
\definecolor{lighter_gray}{gray}{0.5}
\colorlet{light_blue}{blue!20}
\definecolor{dark_green}{rgb}{0.0, 0.6, 0.0}
\definecolor{royal_blue}{rgb}{0.0, 0.22, 0.66}
\definecolor{salmon}{rgb}{1.0, 0.55, 0.41}
\definecolor{gold}{rgb}{0.8, 0.63, 0.21}
\definecolor{navy_blue}{rgb}{0.0, 0.0, 0.5}
\definecolor{crimson}{rgb}{0.79, 0.0, 0.09}
\definecolor{amethyst}{rgb}{0.6, 0.4, 0.8}
\definecolor{alizarin}{rgb}{0.82, 0.1, 0.26}
\definecolor{amaranth}{rgb}{0.9, 0.17, 0.31}
\definecolor{azure}{rgb}{0.0, 0.5, 1.0}
\definecolor{canaryyellow}{rgb}{0.82, 0.41, 0.12}
\definecolor{carrotorange}{rgb}{0.8, 0.33, 0.0}
\definecolor{cadmiumgreen}{rgb}{0.0, 0.42, 0.24}
\definecolor{copper}{rgb}{0.72, 0.45, 0.2}
\definecolor{aqua}{rgb}{0.5, 1.0, 0.83}
\definecolor{awesome}{rgb}{1.0, 0.13, 0.32}
\definecolor{candyapplered}{rgb}{1.0, 0.03, 0.0}
\definecolor{caribbeangreen}{rgb}{0.0, 0.8, 0.6}
\definecolor{aliceblue}{rgb}{0.94, 0.97, 1.0}
\definecolor{babyblue}{rgb}{0.54, 0.81, 0.94}

\newcommand{\bb}{{\boldsymbol b}}

\newcommand{\me}{E\in \mathcal{E}_h}
\newcommand{\bn}{\boldsymbol n} 
\newcommand{\bt}{\boldsymbol t} 
 
\newcommand{\cF}{{\mathcal F}}

\newcommand{\mL}{{\mathrm L}}
\newcommand{\mH}{{\mathrm H}}
\newcommand{\mW}{{\mathrm W}}
\newcommand{\mD}{{\mathrm D}}
\newcommand{\mN}{{\mathrm N}}
\newcommand{\AS}{{\mathrm {AS}}}

\definecolor{light_gray}{gray}{0.75}
\colorlet{light_blue}{blue!20}

\newcommand{\blist}{\begin{list}{}{\itemsep0.0ex\parsep0.1ex\topsep0.2ex\leftmargin1.6em\labelwidth1.3em}}

\newtheorem{lemma}{Lemma}

\newtheorem{remark}[lemma]{Remark}

\definecolor{light_gray}{gray}{0.75}
\definecolor{lighter_gray}{gray}{0.5}
\colorlet{light_blue}{blue!20}
\definecolor{dark_green}{rgb}{0.0, 0.6, 0.0}
\definecolor{royal_blue}{rgb}{0.0, 0.22, 0.66}
\definecolor{salmon}{rgb}{1.0, 0.55, 0.41}
\definecolor{gold}{rgb}{0.8, 0.63, 0.21}
\definecolor{navy_blue}{rgb}{0.0, 0.0, 0.5}
\definecolor{crimson}{rgb}{0.79, 0.0, 0.09}
\definecolor{amethyst}{rgb}{0.6, 0.4, 0.8}
\definecolor{alizarin}{rgb}{0.82, 0.1, 0.26}
\definecolor{amaranth}{rgb}{0.9, 0.17, 0.31}
\definecolor{azure}{rgb}{0.0, 0.5, 1.0}
\definecolor{canaryyellow}{rgb}{0.82, 0.41, 0.12}
\definecolor{carrotorange}{rgb}{0.8, 0.33, 0.0}
\definecolor{cadmiumgreen}{rgb}{0.0, 0.42, 0.24}
\definecolor{copper}{rgb}{0.72, 0.45, 0.2}
\definecolor{aqua}{rgb}{0.5, 1.0, 0.83}
\definecolor{awesome}{rgb}{1.0, 0.13, 0.32}
\definecolor{candyapplered}{rgb}{1.0, 0.03, 0.0}
\definecolor{caribbeangreen}{rgb}{0.0, 0.8, 0.6}

\title{Residual-Based a Posteriori Error Estimators for Algebraic Stabilizations}

\author{Abhinav Jha\footnote{Institute of Applied Analysis and Numerical Simulation, University of Stuttgart,  Pfaffenwaldring 57, 70569,Stuttgart,Germany, \texttt{abhinav.jha@ians.uni-stuttgart.de}}}

\date{}
\begin{document}
\maketitle

\begin{abstract} In this note, we extend the analysis for the residual-based a posteriori error estimators in the energy norm defined for the algebraic flux correction (AFC) schemes \cite{Jha20} to the newly proposed algebraic stabilization schemes \cite{JK21, Kn23}. Numerical simulations on adaptively refined grids are performed in two dimensions showing the higher efficiency of an algebraic stabilization with similar accuracy compared with an AFC scheme.
\end{abstract}

\textbf{keyword} steady-state convection diffusion reaction equations; algebraically stabilized finite element methods; a posteriori estimator; adaptive grid refinement

\section{Introduction}\label{sec:intro}
The convection diffusion reaction equation represents an important mathematical model that defines a wide range of physical phenomena, such as species concentration. The partial differential equation for this model is given by
\begin{alignat}{3}\label{eq:cdr_eqn}
-\varepsilon \Delta u+\bb\cdot \nabla u+cu &=f\qquad &&\mathrm{in}\ \ \Omega\nonumber \\
u&=u_{\mD}\qquad&&\mathrm{on}\ \ \Gamma_{\mD},\nonumber\\
\varepsilon\nabla u\cdot \bn &=g\qquad&&\mathrm{on}\ \ \Gamma_{\mN},
\end{alignat}
where $\varepsilon>0$ is the constant diffusion coeffecient, $\bb$ is the convective field, $c$ is the reaction coeffecient, $f$ is the source/sink term, $u_{\mD}$ is the Dirichlet boundary condition, $g$ is the Neumann boundary condition, $\Omega \subset \mathbb{R}^d$ $d\geq 2$ is a polygonal domain having boundary $\Gamma$ such that $\Gamma=\Gamma_{\mD}\cup \Gamma_{\mN}$ ($\Gamma_{\mD}\cap \Gamma_{\mN}=\emptyset$), $\Gamma_{\mD}$ being the Dirichlet boundary component and $\Gamma_{\mN}$ is the Neumann boundary component.

Whenever convection dominates diffusion, i.e., $L\|\bb\|_{\infty, \Omega}\gg \varepsilon$ (where $L$ is the characteristic length of the problem), we have the presence of layers on the boundary or the interior of the domain. Layers are narrow regions where the solution has a steep gradient. Another important property of Eq.~\eqref{eq:cdr_eqn} is that the solution satisfies maximum principles. Hence, one would like the numerical solution to satisfy its discrete counterpart, the Discrete Maximum Principle (DMP). Standard numerical techniques such as the central finite difference or Galerkin finite element methods fail on both levels, i.e., they do not approximate the layers and fail to satisfy DMP, and hence we get solutions with under and overshoots. 

The linear stabilization schemes such as the Streamline Upwind Petrov Galerkin (SUPG) \cite{BH81} approximate the layers but fail to satisfy the DMP, and hence we have the presence of under and overshoots near the layers \cite{JK07_1}. The \emph{Algebraic Flux Correction} (AFC) scheme \cite{Ku06} belongs to a small class of nonlinear stabilization schemes that satisfy both properties. The numerical analysis of these schemes has been developed in \cite{BJK16}. Recently, a more generalized approach has been presented in \cite{JK21, Kn23} called \emph{Algebraic Stabilizations}. We refer to \cite{BJK22} for a detailed review of these methods.

Another way of resolving the layers is by using adaptive grid refinement which in turn requires the development of a posteriori error estimators. In \cite{Jha20}, a residual-based posteriori estimator was developed in the energy norm for the AFC schemes, which was independent of the choice of limiters. In this work, we extend the analysis of the residual-based a posteriori estimator for the AFC schemes to the algebraically stabilized methods. Numerical studies are two dimensions are performed for the method developed using $\mathbb{P}_1$ finite elements.

\section{Algebraic Stabilization Schemes}\label{sec:algebraic_stab}
In this paper we use standard notions for Sobolev spaces and their norms.  Let $\Omega\subset \mathbb{R}^d$, $d\geq 2$ be a measurable set, then the $\mL^2$-inner product is defined by $\left(\cdot, \cdot\right)$ and its duality pairing by $\langle \cdot,\cdot \rangle$. The norm (semi-norm) on $\mW^{m,p}(\Omega)$ is denoted by $\|\cdot\|_{m,p,\Omega}(|\cdot|_{m,p,\Omega})$ with the convention $\mW^{m,2}(\Omega)=\mH^m(\Omega)$ and $\|\cdot\|_{m,\Omega}=\|\cdot\|_{m,2,\Omega}$ (similarly $|\cdot|_{m,\Omega}=|\cdot|_{m,2,\Omega}$).

Under the assumption
\begin{equation}\label{eq:cdr_cond}
c(x)-\frac{1}{2}\nabla \cdot \bb(x)\geq \sigma>0,
\end{equation}
it is well known that Eq.~\eqref{eq:cdr_eqn} possesses a unique weak solution $u\in  \mH_{\mD}^1(\Omega)$ that satisfies
$$
a(u,v)=\left\langle f,v\right\rangle+\left\langle g,v\right\rangle_{\Gamma_{\mN}}\qquad \forall v\in\mH_{0,\mD}^1(\Omega),
$$
with 
\begin{equation*}
\begin{aligned}
a(u,v)&=\varepsilon \left( \nabla u,\nabla v\right)+\left(\bb\cdot \nabla u,v\right)+(cu,v),\nonumber \\
\mH_{\mD}^1(\Omega)&=\left\lbrace v\in \mH^1(\Omega): v|_{\Gamma_{\mD}}=u_{\mD}\right\rbrace,\nonumber \\
\mH_{0,\mD}^1(\Omega)&=\left\lbrace v\in \mH^1(\Omega): v|_{\Gamma_{\mD}}=0\right\rbrace,
\end{aligned}
\end{equation*}
and $\langle\cdot,\cdot\rangle_{\Gamma_{\mN}}$ is the duality pairing restricted to the Neumann boundary, e.g., see \cite[Sec.~III.1.1]{RST08}.

The algebraic stabilizations schemes for Eq.~\eqref{eq:cdr_eqn} reads as (see \cite{JK21}): Find $u_h\in W_h\left( \subset \mathcal{C}(\overline{\Omega})\cap \mH_{\mD}^1(\Omega)\right)$ such that
\begin{equation}\label{eq:as_scheme}
a_{\AS}(u_h;u_h,v_h)=\left\langle f,v_h\right\rangle+\left\langle g,v_h\right\rangle_{\Gamma_{\mN}}\qquad \forall v_h\in V_h,
\end{equation}
with $a_{\AS}(\cdot;\cdot,\cdot):\mH_{\mD}^1(\Omega)\times\mH_{\mD}^1(\Omega)\times \mH_{0,\mD}^1(\Omega)\rightarrow \mathbb{R}$ and
$$
a_{\AS}(u;v,w)=a(v,w)+b_h(u;v,w),
$$
where $W_h,V_h$ are linear finite-dimensional subspaces of $\mathcal{C}(\overline{\Omega})\cap \mH_{\mD}^1(\Omega)$ and $\mathcal{C}\left(\overline{\Omega}\right)\cap \mH_{0,\mD}^1(\Omega)$, respectively and $b_h(\cdot;\cdot,\cdot)$ is the nonlinear stabilization term.

For the algebraic stabilization scheme proposed in \cite{JK21, Kn23} the stabilization term is given by
\begin{equation}\label{eq:as_stab}
b(u;v,w)=\sum_{i,j=1}^Nb_{ij}(u)v_j w_i\qquad \forall u,v,w\in \mathcal{C}(\overline{\Omega}),
\end{equation}
with
\begin{alignat*}{2}
b_{ij}(u)&=-\max \left\lbrace \left(1-\alpha_{ij}(u)\right)a_{ij}, 0, \left(1-\alpha_{ji}(u)\right)a_{ji} \right\rbrace,\quad  i\neq j\\
 b_{ii}(u)&=-\sum_{i\neq j}b_{ij}(u),
\end{alignat*}
where $N$ is the total number of degree of freedoms, $\alpha_{ij}(u)\in [0,1]$ is the not necessarily 
symmetric limiters, $a_{ij}$ is the stiffness matrix entries corresponding to Eq.~\eqref{eq:cdr_eqn}, and $v_j$ (similarly $w_i$) denote $v(x_j)$.  In this paper we will concentrate on the Symmetric Monotone Upwind-type Algebraically Stabilized (SMUAS) method proposed in \cite{Kn23}. This analysis is also applicable to the Monotone Upwind-type Algebraically Stabilized (MUAS) method proposed in \cite{JK21} as the SMUAS is a generalization of the MUAS method. The SMUAS is an upwind type linearity-preserving method which satisfies the DMP on arbitrary simplicial meshes.  The limiters are computed as follows:
\begin{enumerate}
    \item Compute 
    \begin{alignat*}{2}
       P_i^+  &= \sum_{j\in S_i}|d_{ij}|\left\lbrace (u_i-u_j)^+ + (u_i-u_{ij})^+\right\rbrace,\\
       P_i^-  &=  \sum_{j\in S_i}|d_{ij}|\left\lbrace (u_i-u_j)^- + (u_i-u_{ij})^-\right\rbrace,
     \end{alignat*}
     where $S_i=\{j\in\{1,\ldots,N\}\,:\,a_{ij}\ne 0\}$.
    \item Compute 
       \begin{equation*}
       Q_i^+ = \sum_{j\in S_i}s_{ij}\left\lbrace (u_j-u_i)^+ + (u_{ij}-u_i)^+\right\rbrace, \ \ 
       Q_i^- = \sum_{j\in S_i}s_{ij}\left\lbrace (u_j-u_i)^- + (u_{ij}-u_i)^-\right\rbrace.
    \end{equation*}
    \item Compute
    \begin{equation*}
       R_i^+ = \min\left\{ 1,\frac{Q_i^+}{P_i^+} \right\}, \quad R_i^- =
\min\left\{ 1,\frac{Q_i^-}{P_i^-} \right\},\qquad i=1,\dots, M,
    \end{equation*}
    where $M$ is the number of non-Dirichlet nodes.  If \(P_i^+\) or \(P_i^-\) is zero, one sets \(R_i^+=1\) or \(R_i^-=1\), respectively.  The values of \(R_i^+\) and \(R_i^-\) are set to $1$ for Dirichlet nodes as well. 
    \item Define
    \[ \alpha_{ij} = \begin{cases} 
        R_i^+ & \mbox{ if } u_i>u_j ,\\
        1 & \mbox{ if } u_i=u_j,\\
        R_i^- & \mbox{ if } u_i<u_j,
      \end{cases} \quad  i,j=1,\dots,N.
   \]
\end{enumerate}
In the above definition $s_{ij}=\max\lbrace |a_{ij}|, a_{ji}\rbrace$ and $u_{ij}=u_i + \nabla u_h|_{K_i^j}\cdot (x_i-x_j)$ for all $j\in S_i$. Here $K_i^j$ is a mesh cell containing $x_i$ that is intersected by the half line $\lbrace x_i+\theta (x_i-x_j):\theta >0\rbrace$. In \cite{Kn23}, optimal convergence rates were observed numerically, whereas the algebraic flux correction schemes and the MUAS method did not have the optimal rates on non-Delaunay meshes.

\section{A Posteriori Error Estimators}\label{sec:apost}
The starting point of the analysis is the natural norm of Eq.~\eqref{eq:as_scheme} given by
$$
\|u\|_{\AS}^2=\|u\|_a^2+b_h(u;u,u),
$$
where $\|\cdot \|_a$ is the energy norm given by $ \|u\|_a^2=\varepsilon\|\nabla u\|_{0,\Omega}^2+\sigma \|u\|_{0,\Omega}^2.$

Now, the error in the $\|\cdot \|_{\AS}$ norm can be written as
\begin{alignat}{2}\label{eq:as_norm_calc1}
\|u-u_h\|_{\AS}^2&=a_{\AS}(u_h;u-u_h, u-u_h)\nonumber \\
&=\langle f,u-I_h u\rangle + \langle g, u-I_hu\rangle_{\Gamma_{\mN}} - a(u_h, u-I_hu)\nonumber \\
&\quad +b_h\left(u_h;u-u_h, u-u_h-I_h(u-u_h)\right) + b_h(u_h;u,I_hu-I_h u_h),
\end{alignat}
where we have used Galerkin orthogonality arguments and $I_h$ denotes the Scott-Zhang interpolation operator,\cite[Sec.~4.8]{BS08} (see \cite[Sec.~3.1.1]{Jha20} for a detailed explanation). Now the first two terms in Eq.~\eqref{eq:as_norm_calc1} can be written as follows:
\begin{alignat}{2}\label{eq:as_norm_calc_2}
\|u-u_h\|_{\AS}^2&=\sum_{K\in \mathcal{T}_h}\left( R_K(u_h), u-I_hu\right)_K+\sum_{F\in \mathcal{F}_h}\left\langle R_F(u_h), u-I_hu\right\rangle_F\nonumber \\
&\quad +b_h(u_h;u,I_hu-u_h)+b_h(u_h;u-u_h,u-u_h-I_h(u-u_h)),
\end{alignat}
with $R_K(u_h):=f+\varepsilon \Delta u_u-\bb\cdot \nabla u_u-cu_h|_K$, and
\begin{equation*}
R_F(u_h):=\left\lbrace
\begin{array}{ll}
-\varepsilon \left[\![ \nabla u_h\cdot \bn_F]\!\right]_F& \mathrm{if}\ \ F\in \mathcal{F}_{h,\Omega},\\
g-\varepsilon(\nabla u_h\cdot \bn_F)& \mathrm{if}\ \ F\in \mathcal{F}_{h,\mN},\\
0 & \mathrm{if}\ \ F\in \mathcal{F}_{h,\mD},
\end{array}\right.
\end{equation*}
where $\lbrace\mathcal{T}_h\rbrace_{h>0}$ denotes the triangulation of $\Omega$, $\mathcal{F}_h \left(:=\mathcal{F}_{h,\Omega}\cup \mathcal{F}_{h,\mN}\cup \mathcal{F}_{h,\mD}\right)$ denotes the set of faces (with $\mathcal{F}_{h,\Omega}$, $\mathcal{F}_{h,\mN}$, and $\mathcal{F}_{h,\mD}$ denoting the interior, Neumann, and Dirichlet faces respectively), $\left[\![ \cdot]\!\right]_F$ denotes the jump across the face $F$, and $\bn_F$ denotes the outward pointing normal on face $F$.

The first two terms in Eq.~\eqref{eq:as_norm_calc_2} can be approximated in a standard way using the Cauchy-Schwarz inequality, Young's inequality, interpolation estimates and trace estimates. Using these arguments we get the global upper bound,
\begin{eqnarray}\label{eq:est_all_00}
\lefteqn{\|u-u_h\|_a^2+\frac{C_Y}{C_Y-1}b_h(u_h;u-u_h,u-u_h)}\nonumber\\
& \le &\frac{C^2_Y}{2(C_Y-1)}\sum_{K\in \mathcal{T}_h}\mathrm{min}\left\lbrace \frac{C_I^2}{\sigma},\ \frac{C_I^2h_K^2}{\varepsilon}\right\rbrace\|R_K(u_h)\|_{0,K}^2\nonumber\\
&& + \frac{C^2_Y}{2(C_Y-1)}\sum_{F\in \cF_h} \min \left\{\frac{C_F^2 h_F}{\varepsilon}, \frac{C_F^2}{\sigma^{1/2}\varepsilon^{1/2}} \right\}\| R_F(u_h)\|_{0,F}^2\\
&& + \frac{C_Y}{C_Y-1}b_h(u_h;u,I_hu-u_h) + \frac{C_Y}{C_Y-1}b_h\left(u_h;u-u_h,u-u_h-I_h(u-u_h)\right),\nonumber
\end{eqnarray}
where $C_I$ denotes the interpolation constant \cite[Corollary~4.8.15]{BS08}, $C_F$ is the trace estimate constant \cite[Sec.~3.3]{Ver13}, and $C_Y$ is the generalized Young's inequality constant which will be chosen later.

Using the linearity of $b_h(\cdot;\cdot,\cdot)$ in the second and the third argument, the last two terms on the right-hand side of Eq.~\eqref{eq:est_all_00} can be reduced to 
\begin{eqnarray}\label{eq:d_h_est}
\lefteqn{b_h(u_h;u-u_h,u-u_h-I_h(u-u_h))+b_h(u_h;u,I_h(u-u_h))}\nonumber \\
& = &  b_h(u_h;u-u_h,u-u_h)+b_h(u_h;u_h,I_h(u-u_h)).
\end{eqnarray}
Inserting Eq.~\eqref{eq:d_h_est} in Eq.~\eqref{eq:est_all_00} we notice that the stabilization term on the left and the right cancels. Till here we have followed the same procedure as of \cite{Jha20}. Now we approximate the second term in Eq.~\eqref{eq:d_h_est}. We notice that
\begin{equation}\label{eq:b_h_est}
\begin{aligned}
b_h(u_h;u_h,I_hu-u_h)&=\sum_{i,j=1}^N b_{ij}(u_h)u_{h,j}(I_hu_i-u_{h,i})\nonumber\\
&=\sum_{i\neq j} b_{ij}(u_h)u_{h,j}\left(I_hu_i-u_{h,i}\right)+b_{ii}(u_h)u_{h,i}\left(I_hu_i-u_{h,i}\right)\nonumber\\
&=\sum_{i\neq j}b_{ij}(u_h)\left(u_{h,j}-u_{h,i}\right)\left( I_hu_i-u_{h,i}\right)\nonumber\\
&=\sum_{E\in \mathcal{E}_h}|b_E(u_h)|h_E\left(\nabla u_h\cdot \bt_E,\nabla\left(I_hu-u_h\right)\cdot \bt_E\right)_E,
\end{aligned}
\end{equation}
where $\mathcal{E}_h$ denotes the set of all edges, $h_E$ corresponds to the edge $E$ having end points $x_i$ and $x_j$, $\bt_E$ is the tangential derivative along $E$, and we have used the edge formulation from \cite{BJK22}. 

This term is similar to \cite[Eq.~(35)]{Jha20} and hence we can use the same bound but by replacing $(1-\alpha_E)d_E$ by $|b_E|$ and hence choosing the optimal value of $C_Y=4$,  we get the global upper bound as
\begin{equation}\label{eq:post_upper_bound}
\|u-u_h\|_a^2\leq \eta^2=\eta_1^2+\eta_2^2+\eta_3^2,
\end{equation}
where 
\begin{eqnarray*}
\eta_1^2& := & \sum_{K\in \mathcal{T}_h}\mathrm{min}\left\{ \frac{4 C_I^2}{\sigma},\ \frac{4C_I^2h_K^2}{\varepsilon}\right\}\|R_K(u_h)\|_{0,K}^2, \\ 
\eta_2^2 & :=&   \sum_{F\in \cF_h} \min \left\{\frac{4 C_F^2 h_F}{\varepsilon}, \frac{4 C_F^2}{\sigma^{1/2}\varepsilon^{1/2}} \right\}\| R_F(u_h)\|_{0,F}^2,\\
\eta_3^2 & := & \sum_{\me}\min \Bigg\{\frac{4\kappa_1h_E^2}{
\varepsilon}, \frac{4\kappa_2}{\sigma}\Bigg\} |b_E|^2 h_E^{1-d} \|\nabla u_h \cdot\bt_E\|_{0,E}^2\nonumber.
\end{eqnarray*}
The constants $\kappa_1$ and $\kappa_2$ are given by
\[
\kappa_1=CC_{\mathrm{edge,max}}\left(1+\left(1+C_I\right)^2\right),\quad \kappa_2=CC_{\mathrm{inv}}^2C_{\mathrm{edge,max}}\left( 1+\left(1+C_I\right)^2\right),
\]
where $C$ is a general constant independent of $h$, $C_{\mathrm{inv}}$ is an
inverse inequality constant \cite[Lemma~4.5.3]{BS08}, and $C_{\mathrm{edge,max}}$ is a computable
constant given by \cite[Remark~9]{Jha20}.  

\begin{remark}
The local lower bound for the residual-based estimator was derived in \cite[Theorem~2]{Jha20}. We would like to note that the lower bound was a formal local lower bound as the edge estimates ($\eta_3$) were approximated. We can use the same formal local lower bound here as
$$
|b_E|\leq C\left( \varepsilon+\|\bb\|_{\infty, \Omega}h+\|c\|_{\infty, \Omega}h^2\right)h_E^{d-2},
$$
see \cite[Eq.~(65)]{JK21}.
\end{remark}
\section{Numerical Studies}\label{sec:numres}
This section presents numerical studies for the a posteriori error estimator developed in this paper. We compare the algebraic stabilization schemes introduced in Sec.~\ref{sec:algebraic_stab} to the AFC scheme with the BJK limiter \cite{BJK17}. In \cite{Jha20}, it was shown that the BJK limiter was the most accurate method, although not the most efficient.

The general strategy of solving a posteriori error problem includes
$$
\mbox{\textbf{SOLVE}} \rightarrow \mbox{\textbf{ESTIMATE}} \rightarrow \mbox{\textbf{MARK}} \rightarrow \mbox{\textbf{REFINE}}.
$$
In this paper, we will use red-green refinements \cite{Ver13} to refine the grids, and for the marking of the mesh cells, we use the maximum marking strategy \cite{Ver13}. We use a fixed-point technique referred to as \emph{fixed point rhs} in \cite{JJ19} along with a dynamic damping parameter to solve the system of nonlinear equations. The stopping criteria is set to $10^4$ nonlinear loops or when the residual is less than $10^{-8}\sqrt{\# \mathrm{dofs}}$, where $\# \mathrm{dofs}$ stands for the number of degrees of freedom.

The constants arising in Eq.~\eqref{eq:post_upper_bound} are set to unity during the simulations, like in \cite{Jha20}. For the adaptive loop the stopping criteria is if $\eta\leq 10^{-3}$ or when $\# \mathrm{dofs}\geq 10^5$.  The example considered in this paper is defined on the unit square, i.e., $\Omega=(0,1)^2$, and the initial grid is obtained by joining the line $(0,0)$ and $(1,1)$. We start the simulation by uniformly refining the grid so that $\#\mathrm{dofs}=25$, followed by two steps of uniform refinement and then the adaptive process. All the simulations used $\mathbb{P}_1$ finite elements.  We would like to note that this is the first simulation with SMUAS on adaptively refined grids.

\textbf{Known Boundary Layer:}
In this example we consider Eq.~\eqref{eq:cdr_eqn} with $\varepsilon=10^{-3}$, $\bb=(2,1)^{\top}$, $c=1$, $g=0$, $u_{\mD}=0$, and the right-hand side $f$ such that the exact solution is given by
$$
u(x,y)=y(1-y)\left(x-\frac{e^{\frac{(x-1)}{\varepsilon}}-e^{-1/\varepsilon}}{1-e^{-1/\varepsilon}}\right).
$$
In this example we have a boundary boundary layer on the right boundary of the domain.

\begin{figure}[t!]
\begin{center}
\begin{tikzpicture}[scale=0.5]
\begin{loglogaxis}[
    legend pos=north east, xlabel = \Large{$\#\ \mathrm{dofs}$}, ylabel = \Large{$\|u-u_h\|_a$},
    legend cell align ={left},, title = {\Large{$\varepsilon=10^{-3}$}},
    legend style={nodes={scale=0.75, transform shape}}]
\addplot[color=crimson,  mark=oplus*, line width = 0.5mm, dashdotted,,mark options = {scale= 1.5, solid}]
coordinates{( 81.0 , 0.12474517392538148 )( 289.0 , 0.12820995240551963 )( 344.0 , 0.12629568389407347 )( 423.0 , 0.12190770396508249 )( 602.0 , 0.11290888039327727 )( 956.0 , 0.09532775832941305 )( 1248.0 , 0.08487973243166279 )( 1723.0 , 0.06718453156052927 )( 2753.0 , 0.04904779068192773 )( 3302.0 , 0.04188463154712021 )( 3897.0 , 0.03756292208679059 )( 5740.0 , 0.036702588407757794 )( 8384.0 , 0.02613363846909399 )( 10354.0 , 0.022004909493320093 )( 14033.0 , 0.018198232314792527 )( 18971.0 , 0.014747273266480138 )( 24428.0 , 0.012221653279338265 )( 37102.0 , 0.009770187096671316 )( 56087.0 , 0.007823069442008883 )( 74590.0 , 0.006548618216429832 )( 105154.0 , 0.0053934154885331214 )};
\addlegendentry{\Large{BJK}} 
\addplot[color=blue,  mark=square*, line width = 0.5mm, dashdotted,,mark options = {scale= 1.5, solid}]
coordinates{( 81.0 , 0.12439836933572607 )( 289.0 , 0.12786718783422693 )( 344.0 , 0.12609034578934433 )( 435.0 , 0.12131913221189022 )( 638.0 , 0.11216655850738808 )( 1045.0 , 0.09369127240435632 )( 1728.0 , 0.06725160704851711 )( 2794.0 , 0.048326146946198016 )( 3321.0 , 0.04122558983745831 )( 4040.0 , 0.03729631383285451 )( 5557.0 , 0.03646545276623516 )( 7184.0 , 0.036144899653049865 )( 9934.0 , 0.027487594603276804 )( 12400.0 , 0.022055703472314173 )( 15702.0 , 0.019347984061960952 )( 21100.0 , 0.015391536368860015 )( 28283.0 , 0.011993583155042759 )( 40635.0 , 0.00975422731951895 )( 53769.0 , 0.008324073266984382 )( 72786.0 , 0.00682354278888873 )( 97326.0 , 0.005683873827914207 )( 129078.0 , 0.0049726597244506305 )};
\addlegendentry{\Large{SMUAS}} 
\addplot[color=violet,  line width = 0.5mm, dashdotted,,mark options = {scale= 1.5, solid}]
coordinates{( 81.0 , 0.1111111111111111 )( 289.0 , 0.058823529411764705 )( 344.0 , 0.053916386601719206 )( 423.0 , 0.048621663832631515 )( 602.0 , 0.04075695729696112 )( 956.0 , 0.03234231136765754 )( 1248.0 , 0.02830692585361489 )( 1723.0 , 0.02409114054616049 )( 2753.0 , 0.01905885887736274 )( 3302.0 , 0.017402492911639783 )( 3897.0 , 0.01601897771107376 )( 5740.0 , 0.013199091933711365 )( 8384.0 , 0.010921300708262974 )( 10354.0 , 0.009827564880753334 )( 14033.0 , 0.008441599375370293 )( 18971.0 , 0.007260305376318418 )( 24428.0 , 0.006398173965797621 )( 37102.0 , 0.005191601383050268 )( 56087.0 , 0.004222492574660969 )( 74590.0 , 0.0036615055471347605 )( 105154.0 , 0.0030838063698310396 )};
\addlegendentry{\Large{$\mathcal{O}(h)$}} 
\end{loglogaxis}
\end{tikzpicture}
\begin{tikzpicture}[scale=0.5]
\begin{loglogaxis}[
    legend pos=north east, xlabel = \Large{$\#\ \mathrm{dofs}$}, ylabel = \Large{$\eta$},
    legend cell align ={left},, title = {\Large{$\varepsilon=10^{-3}$}},
    legend style={nodes={scale=0.75, transform shape}},  fill opacity=0.8, draw opacity=1,text opacity=1]
\addplot[color=crimson,  mark=oplus*, line width = 0.5mm, dashdotted,,mark options = {scale= 1.5, solid}] 
coordinates{( 81.0 , 17.6517 )( 289.0 , 23.8946 )( 344.0 , 125.622 )( 423.0 , 102.026 )( 602.0 , 67.4411 )( 956.0 , 39.6547 )( 1248.0 , 30.4467 )( 1723.0 , 17.7852 )( 2753.0 , 8.98763 )( 3302.0 , 5.16045 )( 3897.0 , 1.61075 )( 5740.0 , 0.757503 )( 8384.0 , 0.455964 )( 10354.0 , 0.28236 )( 14033.0 , 0.222574 )( 18971.0 , 0.17643 )( 24428.0 , 0.147087 )( 37102.0 , 0.116938 )( 56087.0 , 0.0929849 )( 74590.0 , 0.0779277 )( 105154.0 , 0.0645431 )};
\addlegendentry{\Large{BJK}} 
\addplot[color=blue,  mark=square*, line width = 0.5mm, dashdotted,,mark options = {scale= 1.5, solid}] 
coordinates{( 81.0 , 18.6963 )( 289.0 , 25.5382 )( 344.0 , 133.745 )( 435.0 , 106.165 )( 638.0 , 68.9469 )( 1045.0 , 39.8021 )( 1728.0 , 21.9197 )( 2794.0 , 11.9516 )( 3321.0 , 7.5248 )( 4040.0 , 3.73505 )( 5557.0 , 2.14759 )( 7184.0 , 1.01182 )( 9934.0 , 0.620944 )( 12400.0 , 0.385436 )( 15702.0 , 0.270725 )( 21100.0 , 0.199776 )( 28283.0 , 0.151077 )( 40635.0 , 0.121155 )( 53769.0 , 0.0995545 )( 72786.0 , 0.081478 )( 97326.0 , 0.0678458 )( 129078.0 , 0.0596085 )};
\addlegendentry{\Large{SMUAS}} 
\addplot[color=violet,  line width = 0.5mm, dashdotted,,mark options = {scale= 1.5, solid}]
coordinates{( 81.0 , 1.5713484026367723 )( 289.0 , 0.8318903308077029 )( 344.0 , 0.7624928516630234 )( 423.0 , 0.6876141641725289 )( 602.0 , 0.5763904177042349 )( 956.0 , 0.4573893537463482 )( 1248.0 , 0.40032038451271784 )( 1723.0 , 0.3407001769341653 )( 2753.0 , 0.26953296707721247 )( 3302.0 , 0.24610841494742636 )( 3897.0 , 0.22654255534352827 )( 5740.0 , 0.18666334823663935 )( 8384.0 , 0.15445051580380387 )( 10354.0 , 0.1389827553946289 )( 14033.0 , 0.11938224324768915 )( 18971.0 , 0.10267622330159805 )( 24428.0 , 0.09048384396853446 )( 37102.0 , 0.07342033086344607 )( 56087.0 , 0.05971506266105231 )( 74590.0 , 0.051781508034622975 )( 105154.0 , 0.04361160791947597 )};
\addlegendentry{\Large{$\mathcal{O}(h)$}} 
\end{loglogaxis}
\end{tikzpicture}
\begin{tikzpicture}[scale=0.5]
\begin{semilogxaxis}[scaled y ticks = false, axis y line*=left,xticklabels={,,},
    legend pos=north east, ylabel = \Large{iterations+rejections},
    legend cell align ={left}, title = {\Large{$\varepsilon=10^{-3}$}},
     legend style={nodes={scale=0.75, transform shape}}]
\addplot[color=crimson,  mark=oplus*, line width = 0.5mm, dashdotted,, mark options = {scale= 1.5, solid}] 
coordinates{( 25.0 , 79.0 )( 81.0 , 173.0 )( 289.0 , 337.0 )( 344.0 , 275.0 )( 423.0 , 493.0 )( 602.0 , 305.0 )( 956.0 , 231.0 )( 1248.0 , 247.0 )( 1723.0 , 238.0 )( 2753.0 , 170.0 )( 3302.0 , 170.0 )( 3897.0 , 131.0 )( 5740.0 , 107.0 )( 8384.0 , 71.0 )( 10354.0 , 67.0 )( 14033.0 , 59.0 )( 18971.0 , 49.0 )( 24428.0 , 44.0 )( 37102.0 , 27.0 )( 56087.0 , 16.0 )( 74590.0 , 12.0 )( 105154.0 , 6.0 )};
\label{HMM86_BJK17_6_H}
\end{semilogxaxis}
\begin{semilogxaxis}[scaled y ticks = false, axis y line*=right,
    legend pos=north east, xlabel = \Large{$\#\ \mathrm{dofs}$}, 
    legend cell align ={left},
     legend style={nodes={scale=0.75, transform shape}}]
\addlegendimage{/pgfplots/refstyle=HMM86_BJK17_6_H}\addlegendentry{\Large{BJK}}
\addplot[color=blue,  mark=square*, line width = 0.5mm, dashdotted,, mark options = {scale= 1.5, solid}] 
coordinates{( 25.0 , 31.0 )( 81.0 , 52.0 )( 289.0 , 64.0 )( 344.0 , 62.0 )( 435.0 , 57.0 )( 638.0 , 55.0 )( 1045.0 , 51.0 )( 1728.0 , 50.0 )( 2794.0 , 48.0 )( 3321.0 , 47.0 )( 4040.0 , 50.0 )( 5557.0 , 45.0 )( 7184.0 , 44.0 )( 9934.0 , 40.0 )( 12400.0 , 40.0 )( 15702.0 , 38.0 )( 21100.0 , 32.0 )( 28283.0 , 28.0 )( 40635.0 , 21.0 )( 53769.0 , 18.0 )( 72786.0 , 11.0 )( 97326.0 , 7.0 )( 129078.0 , 5.0 )};
\addlegendentry{\Large{SMUAS}} 
\end{semilogxaxis}
\end{tikzpicture}
\end{center}
\caption{Error in the energy norm (left), the global upper bound estimator $\eta$ defined in Eq.~\eqref{eq:post_upper_bound} (center), and the number of nonlinear iterations and rejections (right).  In the right figure, the left $y$-axis corresponds to the AFC scheme with the BJK limiter and the right $y$-axis corresponds to the SMUAS method.}\label{fig:error_energy_norm_example_1_ABR17}
\end{figure}
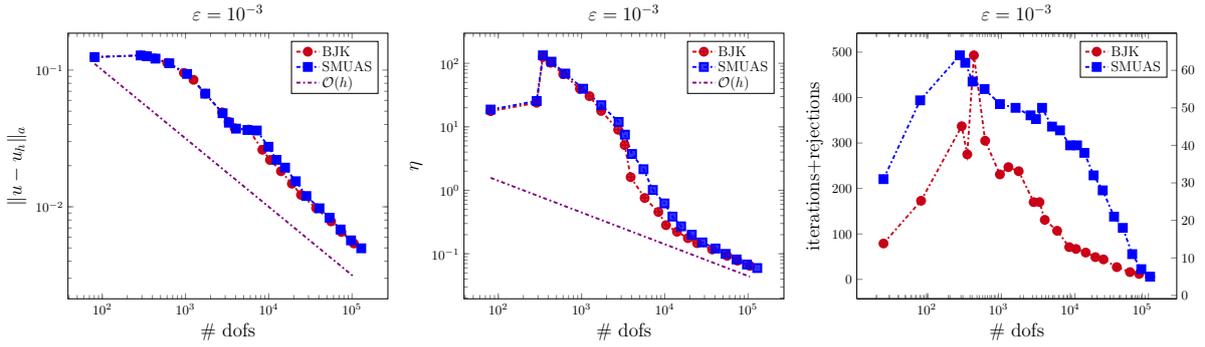

First, in Fig.~\ref{fig:error_energy_norm_example_1_ABR17} we show the error in the energy norm (left) and $\eta$ (center).  We observe that the error decreases at an optimal rate of $\mathcal{O}(h)$ for both the methods, and similarly, so does $\eta$.  We also note the value for the BJK limiter and the SMUAS are comparable.

One of the drawbacks of the AFC scheme with the BJK limiter is the efficiency of the scheme with respect to solving the arising system of equations, see \cite{JJ19, Jha20} for a detailed comparison.  Hence,  we present the result for the number of nonlinear iterations in Fig.~\ref{fig:error_energy_norm_example_1_ABR17} (right). This figure also includes the number of rejections while using dynamic damping. In this figure, the left $y$-axis corresponds to the BJK limiter, whereas the right $y$-axis corresponds to the SMUAS method.  We can observe that SMUAS is very efficient in this regard. In fact it takes one tenth of the number of iterations required by the BJK limiter while providing the same level of accuracy.

\begin{figure}[t!]
\begin{center}
\includegraphics[width=0.3\textwidth]{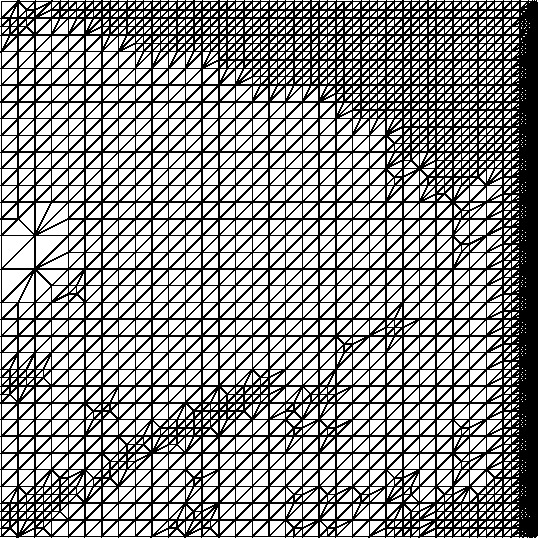}\hspace*{1em}
\includegraphics[width=0.3\textwidth]{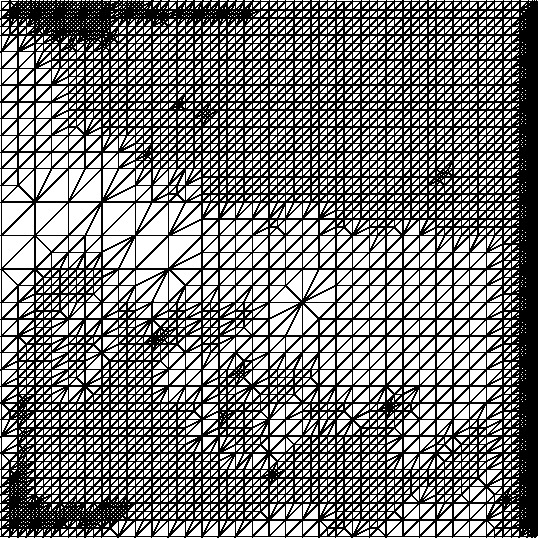}
\end{center}
\caption{Adaptively refined grids with $\#\mathrm{dofs}\approx 25,000$. AFC scheme with BJK limiter (left) and SMUAS method (right).}
\label{fig:conf_grid_example_1_ABR17}
\end{figure}

Finally, in Fig.~\ref{fig:conf_grid_example_1_ABR17}, we present the adaptively refined grids obtained when $\#\mathrm{dofs}\approx 25000$.  We notice that both the methods refine the layer region very much but they are somewhat different in smoother regions.
\section{Summary}\label{sec:summary}
This paper extended the analysis of \cite{Jha20} to the class of new algebraic stabilization schemes, which encompasses the MUAS \cite{JK21} and the SMUAS \cite{Kn23} method. Numerical tests were performed to compare the SMUAS method with the AFC scheme using the BJK limiter. It was observed that both the methods are accurate, and the errors are comparable, but the SMUAS method is much more efficient compared to the BJK limiter. In the future, we would like to perform a deeper numerical study of the different algebraic stabilization schemes on adaptively refined grids.

\bibliographystyle{plain}
\bibliography{Posteriori_MUAS}
\end{document}